\documentclass[smallextended]{article}
\usepackage[english]{babel}
\usepackage[dvips]{graphicx}
\usepackage{amsthm}
\usepackage{exscale}
\usepackage{amsfonts}
\usepackage{amssymb}
\usepackage[intlimits,sumlimits]{amsmath}
\usepackage{amsmath,cases}
\usepackage[latin1]{inputenc}
\usepackage{mathrsfs}
\usepackage{fancyhdr}
\usepackage[misc,geometry]{ifsym}

\theoremstyle{definition}

\newtheorem{thm}{Theorem}

\newcommand{\xib}{\boldsymbol{\xi}}
\newcommand{\thetab}{\boldsymbol{\theta}}
\newcommand{\ud}{\mathrm{d}}

\newcommand{\rone}{\mathbb{R}}

\newcommand{\pp}{\textsf{P}}
\newcommand{\ee}{\textsf{E}}

\newcommand{\samplespace}{\mathbb{X}}
\newcommand{\samplesigmafield}{\mathscr{X}}

\newcommand{\ft}{f(\cdot\ |\ \thetab)}
\newcommand{\TETA}{\thetab \in \Theta}
\newcommand{\ftj}{f(x_j\ |\ \thetab)}

\author{Emanuele Dolera\footnote{Address: Dipartimento di Scienze Fisiche, Informatiche e Matematiche, via Campi 213/b, 41125 Modena, Italy. E-mail: emanuele.dolera@unimore.it} \ and Andrea Bulgarelli\footnote{Address: INAF/IASF-Bologna, via Gobetti 101, I-40129 Bologna, Italy. E-mail: bulgarelli@iasfbo.inaf.it}
}

\title{\textbf{Asymptotic behavior of the joint distribution of a vector of stochastically dependent likelihood ratios}}
\date{}

\begin{document}
\maketitle

\begin{abstract}
This paper provides a generalization of a classical result obtained by Wilks about the asymptotic behavior of the likelihood ratio. The new results deal with the asymptotic behavior of the joint distribution of a vector of likelihood ratios which turn out to be stochastically dependent.
\end{abstract}


\section{Introduction and main results} \label{sect:intro}

The likelihood ratio statistic has been considered, since the pioneering paper \cite{NeyPea}, as a powerful tool to test composite hypotheses, such as $H_0 : \thetab_0 \in \Theta_0$, against the alternative, $H_1 : \thetab_0 \not\in \Theta_0$, where $\Theta_0$ is any proper subset of a generic space of parameters $\Theta$. The basic elements of the analysis are a statistical model, described by the set of probability densities $\{\ft\}_{\TETA}$, and a sample $(\tilde{X}_1, \dots, \tilde{X}_n)$, composed of independent and identically distributed (i.i.d.) random variables (r.v.'s) with density $f(\cdot\ |\ \thetab_0)$. Typically, the true value $\thetab_0$ of the parameter is unknown, and the objective is to draw conclusions on the nature of $\thetab_0$ on the basis of real observations $(x_1, \dots, x_n)$ of the random sample. The quantity
$$
\lambda_n(x_1, \dots, x_n) := \frac{\sup_{\thetab \in \Theta_0} \prod_{j = 1}^{n} \ftj}{\sup_{\thetab \in \Theta} \prod_{j = 1}^{n} \ftj}
$$
corresponds to the aforesaid likelihood ratio statistic, which provides a measure of how well a value $\thetab_0 \in \Theta_0$ explains the observed sample $(x_1, \dots, x_n)$. Indeed, it is reasonable to reject $H_0$ if $\lambda_n$ is too small.
To make these heuristics more effective, one needs the probability distribution of $\lambda_n(\tilde{X}_1, \dots, \tilde{X}_n)$,
under the assumption that $f(\cdot\ |\ \thetab_0)$ defines the common probability law of the $\tilde{X}_i$'s, to specify
the rejection region for a given significance level. A general method to tackle this problem was devised by Wilks in \cite{wilks}, by using asymptotic techniques. Although this procedure is valid only in the presence of large samples and under peculiar assumptions on the model, it has the merit of being very general and of avoiding direct computations. This is particularly useful in presence of complex models, when explicit computations are prohibitive.

This paper is concerned with a generalization of the Wilks theorem, when a vector of likelihood ratios is considered. Suppose, in fact, that $P$ different populations of ``individuals'' are given, and that a sample is drawn from each population. Hence, the data are in the form $(X^{(1)}_1, \dots, X^{(1)}_{n_1})$, $\dots$, $(X^{(P)}_1, \dots, X^{(P)}_{n_P})$
and it is assumed that all these r.v.'s are i.i.d. with common probability law specified by $f(\cdot\ |\ \thetab_0)$. While independence is a very common assumption, given the parameter(s), the identity in distribution is more peculiar, granting that different individuals from distinct populations behave -- from the probabilistic standpoint -- in the same way. For example,   $(X^{(p)}_1, \dots, X^{(p)}_{n_p})$ can be thought of as the outcomes of $n_p$ experiments conducted in the $p$-th laboratory, assuming that all the laboratories work in the same conditions. Otherwise, one can think of $P$ disjoint time intervals and of $(X^{(p)}_1, \dots, X^{(p)}_{n_p})$ as the outcomes of some phenomenon which is observed $n_p$ times in the $p$-th interval, again assuming that the probabilistic description of this phenomenon does not change as time elapses.

The distinctive feature of the present analysis consists in the fact that it is deemed useful to gather in various groups the data coming from the different populations: data from population 1 to population $G$ constitute the first group, data from population 2 to population $G+1$ constitute the second group, and so on, in such a way to form $M$ groups. Of course, $G$ must be an integer in $\{1, \dots, P\}$ and $M = P - G + 1$. This procedure has the following justification. If there are $P$ disjoint time intervals and $(X^{(p)}_1, \dots, X^{(p)}_{n_p})$ represents the outcomes of some phenomenon which is observed $n_p$ times in the $p$-th interval, one can reasonably believe that each time interval is not so representative when considered alone. But these observations can be very significant if groups of $G$ consecutive time intervals are formed, and the data are gathered within each group. In this way, the data belonging to a specific population can be gathered within more than one group, establishing non-null correlations between the $M$ groups. The objective of the analysis remains, in any case, to test the null hypothesis $H_0 : \thetab_0 \in \Theta_0$ against the alternative. Throughout the paper, it will be assumed that $\Theta$ is an open subset of $\rone^d$ containing the origin and that
\begin{equation} \label{eq:H0}
\Theta_0 := \{\mathbf{t} \in \Theta\ |\ t_1 = \dots = t_r = 0\}
\end{equation}
for some integer $r$ in $\{1, \dots, d\}$. Once the data are observed in the form $(x^{(1)}_1, \dots, x^{(1)}_{n_1})$, $\dots$,
$(x^{(P)}_1, \dots, x^{(P)}_{n_P})$ and gathered as above, it is possible to specify the central object of the present analysis, namely the vector of likelihood ratios $(\lambda_{n_1, \dots, n_G}, \dots, \lambda_{n_M, \dots, n_P})$ where
\begin{equation} \label{eq:Lratio}
\lambda_{n_i, \dots, n_{i+G-1}} := \frac{\sup_{\thetab \in \Theta_0} \prod_{p = i}^{i + G - 1} \prod_{j = 1}^{n_p} f(x_j^{(p)}\ |\ \thetab)}{\sup_{\thetab \in \Theta} \prod_{p = i}^{i + G - 1} \prod_{j = 1}^{n_p} f(x_j^{(p)}\ |\ \thetab)}
\end{equation}
for $i = 1, \dots, M$. It is evident from (\ref{eq:H0}) that all the data gathered within the $i$-th group contribute to the evaluation of $\lambda_{n_i, \dots, n_{i+G-1}}$, for all $i = 1, \dots, M$. The achievement of this paper -- encapsulated in Theorem \ref{thm:wilks} -- consists in the evaluation of the limiting probability distribution (as $n_1, \dots, n_P \rightarrow +\infty$) of the vector formed by the various $\lambda$'s, when the data are i.i.d. from $f(\cdot\ |\ \thetab_0)$.

Now, before formalizing the problem in a rigorous mathematical framework, it is worth mentioning the practical problem which originated this kind of procedure. In recent years, a new generation of astronomical observatories (both from ground and from space) has offered great opportunities for discovery in high-energy astrophysics, in particular in the field of short-term variability (a.k.a. transient emission) of $\gamma$-ray astrophysical sources. The AGILE satellite is one of these missions \cite{bulg12}, where the evaluation of the statistical significance of the time variability of these $\gamma$-ray sources is a primary task of the data analysis \cite{cash2}. The scientific operations of the AGILE space mission are focused on the fast detection of these $\gamma$-ray flares, and for this reason an automated analysis systems has been developed \cite{bulg13}, that produce non-independent trials and for this reason the hypothesis of the Wilks theorem \cite{wilks} are not valid.

The mathematical formalization of the situation described above is as follows. Suppose that any observation takes values in a set  $\samplespace$, which is endowed with a $\sigma$-algebra $\samplesigmafield$ and a reference measure $\nu$. Then, a dominated statistical model is given on this space by means of the set of densities $\{\ft\}_{\TETA}$ w.r.t. the dominating measure $\nu$.
Throughout the paper, the space space of the parameters $\Theta$ will be considered as an open subset of $\rone^d$, and the model    will be chosen in conformity to the following regularity conditions.
\begin{enumerate}
\item[i)] $\forall\ x \in \samplespace$, $\thetab \mapsto f(x\ |\ \thetab)$ \emph{is in} $C^2(\Theta; \rone)$;
\item[ii)] \emph{derivatives of first and second order with respect to} $\thetab$ \emph{can be passed under the integral every time it is required};
\item[iii)] \emph{for any fixed} $\thetab_0 \in \Theta$, \emph{there exists a measurable function} $K_0 : \samplespace \rightarrow [0, +\infty)$ \emph{and} $\delta_0 > 0$ \emph{such that}
\begin{eqnarray}
\int_{\samplespace} K_0(x) f(x\ |\ \thetab_0) \nu(\ud x) &<& +\infty \nonumber \\
\sup_{|\thetab - \thetab_0| \leq \delta_0} \Big{|} \frac{\partial^2}{\partial \theta_i \partial \theta_j} \log f(x\ |\ \thetab)\Big{|} &\leq& K_0(x) \ \ \ \ \ (\forall\ x \in \samplespace, i, j = 1, \dots, d); \nonumber
\end{eqnarray}
\item[iv)] \emph{the} Fisher information matrix $\mathrm{I}(\thetab) := (\mathrm{I}_{i,j}(\thetab))_{i,j = 1, \dots, d}$, \emph{defined as}
\begin{equation} \label{eq:fisher}
\mathrm{I}_{i,j}(\thetab) := -\int_{\samplespace} \left(\frac{\partial^2}{\partial \theta_i \partial \theta_j} \log f(x\ |\ \thetab)\right) f(x\ |\ \thetab) \nu(\ud x)\ ,
\end{equation}
\emph{is positive definite at every value of} $\thetab$;
\item[v)] \emph{the validity of the identity} $f(x\ |\ \thetab_1) = f(x\ |\ \thetab_2)$ $\nu$-\emph{almost everywhere in the} $x$-\emph{variable entail} $\thetab_1 = \thetab_2$.
\end{enumerate}
In this framework, the data are modeled by the family of random variables (r.v.'s) of the form $\{X_j^{(p)}\}_{\substack{j = 1, \dots, n_p \\ p = 1, \dots, P}}$, taking values in $\samplespace$. It is implicitly assumed the existence of a suitable base space $(\Omega, \mathscr{F})$ to support all these r.v.'s and, for each $\TETA$, the existence of a probability measure (p.m.) $\pp_{\thetab}$ on $(\Omega, \mathscr{F})$ which makes the observations stochastically independent and identically distributed (i.i.d.) r.v.'s with common distribution generated by $f(x\ |\ \thetab)$, i.e.
$$
\pp_{\thetab}\left[\cap_{p=1}^P \cap_{j=1}^{n_p} \{X_j^{(p)} \in A_j^{(p)}\}\right] = \prod_{p=1}^P \prod_{j=1}^{n_p}
\int_{A_j^{(p)}} f(x\ |\ \thetab) \nu(\ud x) \ \ \ \ \ (A_j^{(p)} \in \samplesigmafield)\ .
$$
The key point hinges on the existence of a distinguished element $\thetab_0 \in \Theta$ as true, but unknown, value of the parameter, which becomes the protagonist of the analysis. In particular, this paper describes a procedure to test
the null hypothesis $H_0 : \thetab_0 \in \Theta_0 \subset \Theta$ against the alternative $H_1 : \thetab_0 \not\in \Theta_0$, with $\Theta_0$ as in (\ref{eq:H0}), based on the observation of the vector $\tilde{\boldsymbol{\Lambda}}_{n_1, \dots, n_P} := (\tilde{\lambda}_{n_1, \dots, n_G}, \dots, \tilde{\lambda}_{n_M, \dots, n_P})$ of likelihood ratios with
\begin{eqnarray}
&& \tilde{\lambda}_{n_i, \dots, n_{i + G - 1}} := \frac{\sup_{\thetab \in \Theta_0} \prod_{p = i}^{i + G - 1} \prod_{j = 1}^{n_p} f(X_j^{(p)}\ |\ \thetab)}{\sup_{\thetab \in \Theta} \prod_{p = i}^{i + G - 1} \prod_{j = 1}^{n_p} f(X_j^{(p)}\ |\ \thetab)} \nonumber \\
&=& \frac{L_{n_i, \dots, n_{i + G - 1}}(\thetab_{n_i, \dots, n_{i + G - 1}}^{\ast}; X_1^{(i)}, \dots X_{n_i}^{(i)}; \dots; X_1^{(i + G - 1)}, \dots X_{n_{i+G-1}}^{(i+G-1)})}{L_{n_i, \dots, n_{i + G - 1}}(\hat{\thetab}_{n_i, \dots, n_{i + G - 1}}; X_1^{(i)}, \dots X_{n_i}^{(i)}; \dots; X_1^{(i + G - 1)}, \dots X_{n_{i+G-1}}^{(i+G-1)})} \nonumber
\end{eqnarray}
for $i = 1, \dots, M$, where $\thetab_{n_i, \dots, n_{i + G - 1}}^{\ast}$ ($\hat{\thetab}_{n_i, \dots, n_{i + G - 1}}$, respectively) stands for the maximum likelihood estimator (MLE) over $\Theta_0$ ($\Theta$, respectively). As recalled in the introduction, the analysis of the case $M = 1$ is well-known: after the introduction, by Neyman and Pearson \cite{NeyPea},
of the concept of likelihood ratio as a useful statistic to test composite hypotheses, the use of this very same statistic became very popular with the discovery, due to Wilks \cite{wilks}, of the limiting distribution of each $\tilde{\lambda}_{n_i, \dots, n_{i + G - 1}}$ when the sample size $\sum_{k = i}^{i + G - 1} n_k$ goes to infinity. The main result of this paper generalizes the aforementioned Wilks theorem by providing the \emph{joint limiting distribution} of the random vector $\tilde{\boldsymbol{\Lambda}}_{n_1, \dots, n_P}$ in the case that $M > 1$, when each sample size $n_i$ goes to infinity. It is worth stressing that, at the level of finitely many data, the various ratios can be stochastically dependent, since  $\tilde{\lambda}_{n_i, \dots, n_{i + G -1}}$ and $\tilde{\lambda}_{n_j, \dots, n_{j + G -1}}$ are formed, in part, by using the same data when $|i - j| < G$. It will be shown in the next two proposition that the correlation between $\tilde{\lambda}_{n_i, \dots, n_{i + G -1}}$ and $\tilde{\lambda}_{n_j, \dots, n_{j + G -1}}$ is maintained or not in the limit according on whether $\lim_{n_1, \dots, n_P \rightarrow +\infty} \frac{\sum_{p = a(i,j)}^{b(i,j)} n_p}{\sqrt{\sum_{q = i}^{i+G-1}\sum_{l = j}^{j+G-1} n_q n_l}}$ exists and it is equal to a strictly positive constant or it is equal to zero, where $a(i,j) := \max\{i, j\}$ and $b(i,j) := \min\{i, j\} + G - 1$. Then, to deduce the exact form of the limiting distribution of $\tilde{\boldsymbol{\Lambda}}_{n_1, \dots, n_P}$, the first step consists in providing a result of asymptotic normality for the vector of MLE's, encapsulated in the following
\begin{thm}\label{thm:cramer}
\emph{Assume that the regularity conditions} i)-v) \emph{are in force and that the MLE's} $\hat{\thetab}_{n_i, \dots, n_{i + G - 1}}$ \emph{actually exist as points of} $\Theta$, \emph{for} $i = 1, \dots, M$. \emph{Moreover, fix a distinguished value} $\thetab_0$ \emph{as true value of the parameter, so that the} $X_j^{(p)}$'s \emph{turn out to be i.i.d. r.v.'s according to}
$\pp_{\thetab_0}$. \emph{Then, the joint distribution of the random vector}
$$
\left(\sqrt{\sum_{k = 1}^G n_k} \cdot (\hat{\thetab}_{n_1, \dots, n_G} - \thetab_0), \dots, \sqrt{\sum_{k = M}^P n_k} \cdot (\hat{\thetab}_{n_M, \dots, n_P} - \thetab_0)\right)
$$
\emph{converges weakly to the} $Md$-\emph{dimensional Gaussian distribution with zero means and covariance matrix}
$$
\left(\begin{array}{c|c|c|c}
\rho_{1,1} \mathrm{I}(\thetab_0)^{-1} & \rho_{1,2} \mathrm{I}(\thetab_0)^{-1} & \ldots & \rho_{1,M} \mathrm{I}(\thetab_0)^{-1} \\
\hline
\rho_{2,1} \mathrm{I}(\thetab_0)^{-1} & \rho_{2,2} \mathrm{I}(\thetab_0)^{-1} & \ldots & \rho_{2,M} \mathrm{I}(\thetab_0)^{-1} \\
\hline
\vdots & \vdots & \ddots & \vdots \\
\hline
\rho_{M,1} \mathrm{I}(\thetab_0)^{-1} & \rho_{M,2} \mathrm{I}(\thetab_0)^{-1} & \ldots & \rho_{M,M} \mathrm{I}(\thetab_0)^{-1}
\end{array} \right)
$$
\emph{where} $\rho_{i,i} := 1$, $\rho_{i,j} := 0$ \emph{if} $|i - j| \geq G$, \emph{and} $\rho_{i,j} := \lim_{n_1, \dots, n_P \rightarrow +\infty} \frac{\sum_{p = a(i,j)}^{b(i,j)} n_p}{\sqrt{\sum_{q = i}^{i+G-1}\sum_{l = j}^{j+G-1} n_q n_l}}$, \emph{with}
$a(i,j) := \max\{i, j\}$ \emph{and} $b(i,j) := \min\{i, j\} + G - 1$.
\end{thm}
After establishing this result, it is possible to state the main achievement of the paper apropos of the asymptotic behavior of the probability distribution, evaluated under $\pp_{\thetab_0}$, of the random vector $\tilde{\boldsymbol{\Lambda}}_{n_1, \dots, n_P}$.
\begin{thm}\label{thm:wilks}
\emph{Assume that} $\Theta_0$ \emph{has the form given by} (\ref{eq:H0}). \emph{Then, under the same assumptions of} Theorem \ref{thm:cramer}, \emph{one has that the probability distribution of the random vector} $(-2 \log[\tilde{\lambda}_{n_1, \dots, n_G}], \dots, -2 \log[\tilde{\lambda}_{n_M, \dots, n_P}])$ \emph{converges weakly to the p.d. of}
$$
\left(\sum_{h = 1}^r \xi_{h; 1}^2, \sum_{h = 1}^r \xi_{h; 2}^2, \dots, \sum_{h = 1}^r \xi_{h; M}^2 \right)
$$
\emph{where the} $\xi$'s \emph{possesses the following distributional properties}
\begin{eqnarray}
\xi_{h; i} &\sim& \mathcal{N}(0, 1) \ \ \ \ \ \forall\ h = 1, \dots, r, \forall i = 1, \dots, M \nonumber \\
\textsf{Cov}(\xi_{h; i}, \xi_{l; j}) &=& 0 \ \ \ \ \ \text{if}\ h \neq l, \forall i, j = 1, \dots, M \nonumber \\
\textsf{Cov}(\xi_{h; i}, \xi_{h; j}) &=& 0 \ \ \ \ \ \text{if}\ |i - j| \geq G, \forall\ h = 1, \dots, r  \nonumber \\
\textsf{Cov}(\xi_{h; i}, \xi_{h; j}) &=& \rho_{i,j} \ \ \ \ \ \text{if}\ |i - j| < G, \forall\ h = 1, \dots, r  \nonumber
\end{eqnarray}
$\mathcal{N}(0, 1)$ \emph{standing for the standard 1-dimensional normal distribution}.
\end{thm}
The importance of this last result is evident, since the limiting distribution turns out to be independent of the specific value of $\thetab_0$, which is fixed but always unknown. Then, it is possible to express the exact (limiting) value of probabilities of the form $\pp_{\thetab_0}[\tilde{\lambda}_{n_1, \dots, n_G} > z_1, \dots, \tilde{\lambda}_{n_M, \dots, n_P} > z_M]$ for every $(z_1, \dots, z_M) \in [0, 1]^M$, which provide the probability of first-type error.

\section{Proofs} \label{sect:proofs}

Gathered here are the proofs of Theorems \ref{thm:cramer} and \ref{thm:wilks}.

\subsection{Proof of Theorem \ref{thm:cramer}} \label{sect:proof1}

Start by noting that the existence of the MLE's as points of $\Theta$, which is an open set, entails that $\hat{\thetab}_{n_i, \dots, n_{i + G - 1}}$ can be expressed as a root of the likelihood equation
\begin{eqnarray}
\ell^{'}_{n_i, \dots, n_{i + G - 1}}(\thetab) &:=& \nabla_{\thetab} \log [L_{n_i, \dots, n_{i + G - 1}}(\thetab; X_1^{(i)}, \dots X_{n_i}^{(i)}; \dots; X_1^{(i + G - 1)}, \dots X_{n_{i+G-1}}^{(i+G-1)})] \nonumber \\
&=& \sum_{p = i}^{i + G - 1} \sum_{j = 1}^{n_p} \nabla_{\thetab} \log[f(X_j^{(p)}\ |\ \thetab)] = \mathbf{0} \nonumber
\end{eqnarray}
for every $i = 1, \dots, M$. Under the assumptions of the theorem, it is well-known that these estimators are strongly consistent, that is $\hat{\thetab}_{n_i, \dots, n_{i + G - 1}} \rightarrow \thetab_0$, $\pp_{\thetab_0}$-a.s.. See, for example, the beginning of the proof of Theorem 18 in \cite{ferg}. Now, expand $\ell^{'}_{n_i, \dots, n_{i + G - 1}}$ as
$$
\ell^{'}_{n_i, \dots, n_{i + G - 1}}(\thetab) = \ell^{'}_{n_i, \dots, n_{i + G - 1}}(\thetab_0) + \int_0^1 \left\{\sum_{p = i}^{i + G - 1} \sum_{j = 1}^{n_p} \mathrm{M}(X_j^{(p)}; \thetab_0 + u(\thetab - \thetab_0))\right\} \ud u \cdot (\thetab - \thetab_0)
$$
where $\mathrm{M}(x; \mathbf{t})$ is the $d \times d$ matrix given by $\left(\frac{\partial^2}{\partial t_k \partial t_h} \log f(x\ |\ \mathbf{t})\right)_{k, h = 1, \dots, d}$. Now, let $\thetab = \hat{\thetab}_{n_i, \dots, n_{i + G - 1}}$, where $\hat{\thetab}_{n_i, \dots, n_{i + G - 1}}$ can be any root of the likelihood equation, and divide by $\sqrt{\sum_{k = i}^{i + G - 1} n_k}$ to obtain
$$
\frac{1}{\sqrt{\sum_{k = i}^{i + G - 1} n_k}} \ell^{'}_{n_i, \dots, n_{i + G - 1}}(\thetab_0) = \sqrt{\sum_{k = i}^{i + G - 1} n_k} \cdot \mathrm{B}_{n_i, \dots, n_{i + G - 1}} (\hat{\thetab}_{n_i, \dots, n_{i + G - 1}} - \thetab_0)
$$
with
$$
\mathrm{B}_{n_i, \dots, n_{i + G - 1}} := - \int_0^1 \frac{1}{\sum_{k = i}^{i + G - 1} n_k} \left\{\sum_{p = i}^{i + G - 1} \sum_{j = 1}^{n_p} \mathrm{M}(X_j^{(p)}; \thetab_0 + u(\hat{\thetab}_{n_i, \dots, n_{i + G - 1}} - \thetab_0))\right\} \ud u\ .
$$
It is well-known that $\mathrm{B}_{n_i, \dots, n_{i + G - 1}} \rightarrow \mathrm{I}(\thetab_0)$, $\pp_{\thetab_0}$-a.s., as shown, for example, in the final part of the proof of Theorem 18 in \cite{ferg}. Therefore, the original problem is traced back to the determination of the limiting distribution of the $Md$-dimensional random vector
$$
\mathbf{V}_{n_1, \dots, n_P} := \left(\frac{1}{\sqrt{\sum_{k = 1}^G n_k}} \ell^{'}_{n_1, \dots, n_G}(\thetab_0), \dots, \frac{1}{\sqrt{\sum_{k = M}^P n_k}} \ell^{'}_{n_M, \dots, n_P}(\thetab_0)\right)
$$
where, by definition,
$$
\frac{1}{\sqrt{\sum_{k = i}^{i+G-1} n_k}} \ell^{'}_{n_i, \dots, n_{i+G-1}}(\thetab_0) =
\sqrt{\sum_{k = i}^{i+G-1} n_k} \left(\frac{1}{\sum_{k = i}^{i+G-1} n_k} \sum_{p = i}^{i+G-1} \sum_{j = 1}^{n_p} \mathbf{\Psi}(X_j^{(p)}; \thetab_0)\right)
$$
for $i = 1, \dots, M$, with $\mathbf{\Psi}(x; \mathbf{t}) := \nabla_{\mathbf{t}} \log[f(X_j^{(p)}\ |\ \mathbf{t})]$. It is worth noticing, at this stage, that the random vectors $\{\mathbf{\Psi}(X_j^{(p)}; \thetab_0)\}_{\substack{j = 1, \dots, n_p \\ p = 1, \dots, P}}$ are i.i.d. under $\pp_{\thetab_0}$, and it follows from hypothesis ii) on the model that
\begin{eqnarray}
\ee_{\thetab_0}[\mathbf{\Psi}(X_j^{(p)}; \thetab_0)] &=& \mathbf{0} \label{eq:mediaPsi} \\
\textsf{Cov}_{\thetab_0}(\Psi^{(k)}(X_j^{(p)}; \thetab_0), \Psi^{(h)}(X_j^{(p)}; \thetab_0)) &=& \mathrm{I}_{k, h}(\thetab_0) \label{eq:varianzaPsi}
\end{eqnarray}
where $\Psi^{(k)}(X_j^{(p)}; \thetab_0)$ denotes the $k^{th}$ coordinate of $\mathbf{\Psi}(X_j^{(p)}; \thetab_0)$. For notational convenience, define $\mathbf{S}_p := \sum_{j = 1}^{n_p} \mathbf{\Psi}(X_j^{(p)}; \thetab_0)$ for $p = 1, \dots, P$, and note that they are independent $d$-dimensional random vector, under $\pp_{\thetab_0}$. Then, the characteristic function of the random vector $\mathbf{V}_{n_1, \dots, n_P}$ is given by
\begin{eqnarray}
\Phi_{n_1, \dots, n_P}(\xib_1, \dots, \xib_M) &=& \ee_{\thetab_0}\left[\exp\left\{\sum_{m=1}^M \frac{i\xib_m \bullet \sum_{p=m}^{m+G-1} \mathbf{S}_p}{\sqrt{\sum_{k = m}^{m+G-1} n_k}} \right\}\right] \nonumber \\
&=& \ee_{\thetab_0}\left[\exp\left\{\sum_{p=1}^P \mathbf{S}_p \bullet \sum_{\substack{m = 1, \dots, M \\ m \leq p \leq m+G-1}}
\frac{i\xib_m}{\sqrt{\sum_{k = m}^{m+G-1} n_k}}\right\}\right] \nonumber \\
&=& \prod_{p=1}^P \varphi^{n_p}\left( \sum_{\substack{m = 1, \dots, M \\ m \leq p \leq m+G-1}}
\frac{\xib_m}{\sqrt{\sum_{k = m}^{m+G-1} n_k}}\right) \nonumber
\end{eqnarray}
where $\bullet$ stands for the standard scalar product in $\rone^d$ and $\varphi(\xib) := \ee_{\thetab_0}[\exp\{i \xib \bullet \mathbf{\Psi}(X_j^{(p)}; \thetab_0)\}]$, with $\xib \in \rone^d$. As in the standard proof of the multi-dimensional central limit theorem (see, for example, Proposition 5.9 and Lemma 5.10 in \cite{ka}), one exploits (\ref{eq:mediaPsi})-(\ref{eq:varianzaPsi}) to deduce
\begin{gather}
\Phi_{n_1, \dots, n_P}(\xib_1, \dots, \xib_M) \nonumber \\
= \prod_{p=1}^P \left[1 - \frac{1}{2}\ ^t\boldsymbol{\Xi}^{(p)}_{n_1, \dots, n_P}(\xib_1, \dots, \xib_M) \mathrm{I}(\theta_0) \boldsymbol{\Xi}^{(p)}_{n_1, \dots, n_P}(\xib_1, \dots, \xib_M)
+ o(\frac{1}{n_p})\right]^{n_p} \nonumber
\end{gather}
with
$$
\boldsymbol{\Xi}^{(p)}_{n_1, \dots, n_P}(\xib_1, \dots, \xib_M) := \sum_{\substack{m = 1, \dots, M \\ m \leq p \leq m+G-1}}\frac{\xib_m}{\sqrt{\sum_{k = m}^{m+G-1} n_k}}\ .
$$
Here, $o(\frac{1}{n_p})$ is complex-valued so, by taking the principal branch of the complex logarithm, one has
\begin{gather}
\text{Log}[\Phi_{n_1, \dots, n_P}(\xib_1, \dots, \xib_M)]  \nonumber \\
= - \frac{1}{2} \sum_{p=1}^P n_p\ ^t\boldsymbol{\Xi}^{(p)}_{n_1, \dots, n_P}(\xib_1, \dots, \xib_M) \mathrm{I}(\theta_0) \boldsymbol{\Xi}^{(p)}_{n_1, \dots, n_P}(\xib_1, \dots, \xib_M) + o(1) \nonumber \ .
\end{gather}
The former term in the RHS above is evidently a quadratic form in the $\xib$-variables, which can be written as
\begin{eqnarray}
&& \sum_{p=1}^P n_p \sum_{\substack{m = 1, \dots, M \\ m \leq p \leq m+G-1}} \sum_{\substack{l = 1, \dots, M \\ l \leq p \leq l+G-1}}\left(\frac{n_p}{\sqrt{\sum_{k = m}^{m+G-1} \sum_{h = l}^{l+G-1} n_k n_h}}\right)\ ^t\xib_m \mathrm{I}(\theta_0)\xib_l \nonumber \\
&=& \sum_{\substack{m, l = 1, \dots, M \\ |m-l| < G}} \left(\frac{\sum_{p = a(l,m)}^{b(l, m)}n_p}{\sqrt{\sum_{k = m}^{m+G-1} \sum_{h = l}^{l+G-1} n_k n_h}}\right)\ ^t\xib_m \mathrm{I}(\theta_0)\xib_l \nonumber
\end{eqnarray}
where $a(l,m) := \max\{l, m\}$ and $b(l,m) := \min\{l, m\} + G - 1$. At this stage, taking the limit as $n_1, \dots, n_P \rightarrow +\infty$ of the above expression, one gets
$$
\lim_{n_1, \dots, n_P \rightarrow +\infty} \text{Log}[\Phi_{n_1, \dots, n_P}(\xib_1, \dots, \xib_M)] = - \frac{1}{2} \sum_{\substack{m, l = 1, \dots, M \\ |m-l| < G}} \rho_{l, m} \ ^t\xib_m \mathrm{I}(\theta_0)\xib_l
$$
for every fixed $\xib_1, \dots, \xib_M \in \rone^d$. This fact, in view of the L\'{e}vy continuity theorem, amounts to proving that the probability distribution of $\mathbf{V}_{n_1, \dots, n_P}$, evaluated under $\pp_{\thetab_0}$, converges weakly to the $Md$-dimensional normal distribution with zero means and covariance matrix given by
\begin{equation}\label{eq:covrho}
\left(\begin{array}{c|c|c|c}
\rho_{1,1} \mathrm{I}(\thetab_0) & \rho_{1,2} \mathrm{I}(\thetab_0) & \ldots & \rho_{1,M} \mathrm{I}(\thetab_0) \\
\hline
\rho_{2,1} \mathrm{I}(\thetab_0) & \rho_{2,2} \mathrm{I}(\thetab_0) & \ldots & \rho_{2,M} \mathrm{I}(\thetab_0) \\
\hline
\vdots & \vdots & \ddots & \vdots \\
\hline
\rho_{M,1} \mathrm{I}(\thetab_0) & \rho_{M,2} \mathrm{I}(\thetab_0) & \ldots & \rho_{M,M} \mathrm{I}(\thetab_0)
\end{array} \right) \ .
\end{equation}
Therefore, upon observing that
\begin{gather}
\mathbf{V}_{n_1, \dots, n_P} = \nonumber \\
\left(\sqrt{\sum_{k = 1}^G n_k} \cdot \mathrm{B}_{n_1, \dots, n_G} (\hat{\thetab}_{n_1, \dots, n_G} - \thetab_0), \dots, \sqrt{\sum_{k = M}^P n_k} \cdot \mathrm{B}_{n_M, \dots, n_P} (\hat{\thetab}_{n_M, \dots, n_P} - \thetab_0)\right)\ , \nonumber
\end{gather}
and that $\mathrm{B}_{n_i, \dots, n_{i + G - 1}} \rightarrow \mathrm{I}(\thetab_0)$, $\pp_{\thetab_0}$-a.s., the desired conclusion now follows, via an obvious application of the Slutsky theorem, from the above achievement on the limiting distribution of $\mathbf{V}_{n_1, \dots, n_P}$, thanks to the elementary property of normal distributions.

\subsection{Proof of Theorem \ref{thm:wilks}} \label{sect:proof2}

The argumentation developed at the beginning of the proof of Theorem 22 in \cite{ferg} shows that the limiting distribution of $(-2 \log[\tilde{\lambda}_{n_1, \dots, n_G}], \dots, -2 \log[\tilde{\lambda}_{n_M, \dots, n_P}])$, evaluated under $\pp_{\thetab_0}$, is the same as the limiting distribution of the $M$-dimensional random vector $\mathbf{W}_{n_1, \dots, n_P}$
whose components are given by
$$
\frac{1}{\sqrt{\sum_{k = i}^{i+G-1} n_k}}\ ^t\ell^{'}_{n_i, \dots, n_{i+G-1}}(\thetab_{n_i, \dots, n_{i + G - 1}}^{\ast})
\mathrm{I}(\thetab_0)^{-1} \frac{1}{\sqrt{\sum_{k = i}^{i+G-1} n_k}} \ell^{'}_{n_i, \dots, n_{i+G-1}}(\thetab_{n_i, \dots, n_{i + G - 1}}^{\ast})
$$
for $i = 1, \dots, M$, $\thetab_{n_i, \dots, n_{i + G - 1}}^{\ast}$ denoting the MLE over $\Theta_0$ based on the observations
$\{X_j^{(p)}\}_{\substack{j = 1, \dots, n_p \\ p = i, \dots, i+G-1}}$. But, exactly as in the central part of the above-mentioned proof from \cite{ferg}, the limiting distribution of the random vector
$$
\left(\frac{1}{\sqrt{\sum_{k = 1}^G n_k}} \ell^{'}_{n_1, \dots, n_G}(\thetab_{n_1, \dots, n_G}^{\ast}), \dots,
\frac{1}{\sqrt{\sum_{k = M}^P n_k}} \ell^{'}_{n_M, \dots, n_P}(\thetab_{n_M, \dots, n_P}^{\ast})\right)
$$
turns out to be the same as the limiting distribution of the random vector with components given by
$$
[\mathrm{Id}_{d\times d} - \mathrm{I}(\thetab_0)\mathrm{H}(\thetab_0)] \frac{1}{\sqrt{\sum_{k = i}^{i+G-1} n_k}} \ell^{'}_{n_i, \dots, n_{i+G-1}}(\thetab_0)
$$
for $i = 1, \dots, M$, where $\mathrm{H}(\thetab_0)$ is a $d\times d$ matrix defined as follows. Partition $\mathrm{I}(\thetab_0)$ into four matrices in such a way that
$$
\mathrm{I}(\thetab_0) = \left(\begin{array}{c|c}
\mathrm{G}_1(\thetab_0) & \mathrm{G}_2(\thetab_0) \\
\hline
^t\mathrm{G}_2(\thetab_0) & \mathrm{G}_3(\thetab_0)
\end{array} \right)
$$
holds with $\mathrm{G}_1(\thetab_0)$ of dimension $r\times r$, $\mathrm{G}_2(\thetab_0)$ of dimension $r\times (d-r)$ and
$\mathrm{G}_3(\thetab_0)$ of dimension $(d-r) \times (d-r)$, and let
$$
\mathrm{H}(\thetab_0) := \left(\begin{array}{c|c}
\mathbf{0} & \mathbf{0} \\
\hline
\mathbf{0} & \mathrm{G}_3(\thetab_0)^{-1}
\end{array} \right)\ .
$$
At this stage, it is possible to exploit the form of the limiting distribution of the random vector
$\mathbf{V}_{n_1, \dots, n_P}$ deduced in the proof of Theorem \ref{thm:cramer}. Indeed, an application of the continuous mapping theorem entails that the limiting distribution of $(-2 \log[\tilde{\lambda}_{n_1, \dots, n_G}], \dots, -2 \log[\tilde{\lambda}_{n_M, \dots, n_P}])$ coincides with the probability law of the random vector with component given by
$$
^t\mathbf{Y}_i\ ^t[\mathrm{Id}_{d\times d} - \mathrm{I}(\thetab_0)\mathrm{H}(\thetab_0)]\ \mathrm{I}(\thetab_0)^{-1}\ [\mathrm{Id}_{d\times d} - \mathrm{I}(\thetab_0)\mathrm{H}(\thetab_0)]\ \mathbf{Y}_i
$$
for $i = 1, \dots, M$, where $(^t\mathbf{Y}_1, \dots, ^t\mathbf{Y}_M)$ is any $Md$-dimensional random vector having normal distribution with zero means and covariance matrix (\ref{eq:covrho}). Elementary linear algebra shows that
$$
^t[\mathrm{Id}_{d\times d} - \mathrm{I}(\thetab_0)\mathrm{H}(\thetab_0)] \mathrm{I}(\thetab_0)^{-1} [\mathrm{Id}_{d\times d} - \mathrm{I}(\thetab_0)\mathrm{H}(\thetab_0)] = \mathrm{I}(\thetab_0)^{-1} - \mathrm{H}(\thetab_0)
$$
because $\mathrm{H}(\thetab_0)\mathrm{I}(\thetab_0)\mathrm{H}(\thetab_0) = \mathrm{H}(\thetab_0)$. To conclude, introduce the the  $Md$-dimensional random vector $(^t\mathbf{Z}_1, \dots, ^t\mathbf{Z}_M)$, defined by putting $\mathbf{Z}_i := \mathrm{I}(\thetab_0)^{-1/2} \mathbf{Y}_i$, whose distribution is normal with zero means and covariance matrix equal to
$$
\left(\begin{array}{c|c|c|c}
\rho_{1,1} \mathrm{Id}_{d \times d} & \rho_{1,2} \mathrm{Id}_{d \times d} & \ldots & \rho_{1,M} \mathrm{Id}_{d \times d} \\
\hline
\rho_{2,1} \mathrm{Id}_{d \times d} & \rho_{2,2} \mathrm{Id}_{d \times d} & \ldots & \rho_{2,M} \mathrm{Id}_{d \times d} \\
\hline
\vdots & \vdots & \ddots & \vdots \\
\hline
\rho_{M,1} \mathrm{Id}_{d \times d} & \rho_{M,2} \mathrm{Id}_{d \times d} & \ldots & \rho_{M,M} \mathrm{Id}_{d \times d}
\end{array} \right) \ .
$$
Upon noticing that
$$
\mathrm{I}(\thetab_0)^{1/2}\ [\mathrm{I}(\thetab_0)^{-1} - \mathrm{H}(\thetab_0)]\ \mathrm{I}(\thetab_0)^{1/2} = \mathrm{P}_{d,r} := \left(\begin{array}{c|c}
\mathrm{Id}_{r \times r} & \mathbf{0} \\
\hline
\mathbf{0} & \mathbf{0}
\end{array} \right)
$$
holds, the theorem is completely proved by putting
$$
(\xi_{1;i}, \dots, \xi_{r;i}) := \mathrm{P}_{d,r} \mathbf{Z}_i
$$
for $i = 1, \dots, M$.

\end{document}